\documentclass{article}
\usepackage[skip=-1.5pt]{caption}
\usepackage{tikz}
\usepackage{algorithm}
\usepackage[margin=2.25cm]{geometry}
\begin{document}

\title{\huge On the lattice formulation of the union-closed sets conjecture}
\author{\Large Christopher Bouchard}
\date{}
\maketitle

\smallskip

\abstract{\noindent The union-closed sets conjecture, also known as Frankl's conjecture, is a well-studied problem with various formulations. In terms of lattices, the conjecture states that every finite lattice $L$ with more than one element contains a join-irreducible element that is less than or equal to at most half of the elements in $L$. In this work, we obtain several necessary conditions for any counterexample $\tilde{L}$ of minimum size.}

\bigskip

\section*{1. Introduction}

Let $\mathcal{A}$ be a finite family of distinct finite sets with at least one nonempty member set. The union-closed sets conjecture, also known as Frankl's conjecture, states that if $X,Y \in \mathcal{A}$ implies that $X \cup Y \in \mathcal{A}$, then there exists an element of $\bigcup_{A \in \mathcal{A}}A$ that is in at least $\frac{|\mathcal{A}|}{2}$ member sets of $\mathcal{A}$. With respect to this standard formulation of the conjecture, minimal counterexamples have been studied. For example, in [7] it was shown that for any counterexample $\tilde{\mathcal{A}}$ of minimum size $|\tilde{\mathcal{A}}|$, there must exist at least three elements of $\bigcup_{A \in \tilde{\mathcal{A}}}A$ that are each in exactly $\frac{|\tilde{\mathcal{A}}|-1}{2}$ member sets of $\tilde{\mathcal{A}}$. Another result, from [6] and [10], is that $|\tilde{\mathcal{A}}| \geq 4|\bigcup_{A \in \mathcal{A}^*}A|-1$, where $\mathcal{A}^*$ is a counterexample, in this case, with minimum $|\bigcup_{A \in \mathcal{A}^*}A|$. In the present work, we consider the lattice formulation of the union-closed sets conjecture (found, for instance, in [8] and [11]), obtaining several necessary conditions for any counterexample of minimum size.

To begin, we recall some terminology. A \textit{poset} (or \textit{partially ordered set}) $\langle P; \leq_P \rangle$ is a set $P$ equipped with a binary relation $\leq_P$ on $P$ that is reflexive, antisymmetric, and transitive. For $x, y \in P$, if $x \leq_P y$, then we state that $x$ is less than or equal to $y$ (and $y$ is greater than or equal to $x$) in $P$. Similarly, if $x \leq_P y$ and $x \neq y$ (denoted by $x <_P y$), then we state that $x$ is less than $y$ (and $y$ is greater than $x$) in $P$. Two elements $x, y \in P$ are \textit{comparable} in $P$ if $x \leq_P y$ or $y \leq_P x$. Otherwise, $x$ and $y$ are \textit{incomparable} in $P$, denoted by $x \parallel_P y$ (or $y \parallel_P x$). We define a \textit{subposet} $\langle S; \leq_S \rangle$ of $P$ to be a poset such that $S \subseteq P$ and, for all $x$ and $y$ in $S$, $x \leq_S y$ if and only if $x \leq_P y$. A poset is \textit{linearly ordered} if every two of its elements are comparable, and a \textit{chain} $\langle C; \leq_C \rangle$ in $P$ is a linearly ordered subposet of $P$. We use the notation $({\uparrow}x)_P = \{y \in P \ | \ x \leq_P y \}$, $({\downarrow}x)_P = \{y \in P \ | \ y \leq_P x \}$, and $({\parallel}x)_P=\{y \in P \ | \ x \parallel_P y \}$. Two elements $x, y \in P$ have a \textit{join} $\sup_P\{x,y\}$ in $P$ if $\sup_P\{x,y\}$ is an element of $P$ such that both $x \leq_P \sup_P\{x,y\}$ and $y \leq_P \sup_P\{x,y\}$, and $x \leq_P z$ and $y \leq_P z$ together imply that $\sup_P\{x,y\} \leq_P z$ for all $z \in P$. Similarly, $x$ and $y$ have a \textit{meet} $\inf_P\{x,y\}$ in $P$ if $\inf_P\{x,y\}$ is an element of $P$ such that $\inf_P\{x,y\} \leq_P x$ and $\inf_P\{x,y\} \leq_P y$, and $z \leq_P x$ and $z \leq_P y$ imply that $z \leq_P \inf_P\{x,y\}$. An element $x$ \textit{upper covers} an element $y$ (and $y$ \textit{lower covers} $x$) in $P$ if $y <_P x$, and for all $z \in P$, $y <_P z \leq_P x$ implies that $z = x$. An element is \textit{join-irreducible} in $P$ if it upper covers exactly one element in $P$, \textit{meet-irreducible} if it lower covers exactly one element, and \textit{doubly irreducible} if it is both join and meet-irreducible. On the other hand, an element is \textit{join-reducible} in $P$ if it upper covers more than one element in $P$, \textit{meet-reducible} if it lower covers more than one element, and \textit{doubly reducible} if it is both join and meet-reducible.

A poset $\langle L; \leq_L \rangle$ is a \textit{lattice} if every two elements $x, y \in L$ have both a join and meet in $L$. When $L$ is finite, its greatest element is denoted by $1_L$ and least element by $0_L$. An element that upper covers $0_L$ in $L$ is an \textit{atom}, and an element that lower covers $1_L$ is a \textit{dual atom}. For more information on lattices, see [2]. The lattice formulation of the union-closed sets conjecture (phrased, as usual, according to the intersection-closed dual) is stated as follows:

\medskip

\noindent \textbf{Conjecture 1.1.} \textit{Any finite lattice} $L$ \textit{with more than one element contains a join-irreducible element} $j$ \textit{such that} $|({\uparrow}j)_L| \leq \frac{|L|}{2}$\textit{.}

\medskip

In the next section, we characterize any counterexample $\tilde{L}$ to Conjecture 1.1 of minimum size $|\tilde{L}|$. We conclude the current section by considering some pertinent statements regarding lattices in general.

\medskip

\noindent \textbf{Lemma 1.2} (Agalave, Shewale, and Kharat [1])\textbf{.} \textit{An element} $x \in L$ \textit{is join or meet-irreducible in a finite lattice} $L$ \textit{with} $|L| > 1$ \textit{if and only if the subposet} $L \setminus \{x\}$ \textit{of} $L$ \textit{is a lattice.}

\medskip

\noindent \textit{Proof.} See Lemma 2.1 of [1]. We include here a proof using the present notation. First, we show that if $x$ is join-irreducible in $L$, upper covering a unique element $z$ in $L$, then the subposet $L \setminus \{x\}$ of $L$ is a lattice. Consider any $y_1,y_2 \in L \setminus \{x\}$. If $\inf_L\{y_1,y_2\}=x$, then $\inf_{L \setminus \{x\}}\{y_1,y_2\}=z$. (Otherwise, there exists an element $l$ less than or equal to $y_1$ and $y_2$ in $L \setminus \{x\}$ such that $z <_{L \setminus \{x\}} l$ or $z \parallel_{L \setminus \{x\}} l$ (with $l <_L x$ because $\inf_L\{y_1,y_2\}=x$ and $l \neq x$). However, $z <_{L \setminus \{x\}} l$ does not hold (because $z$ lower covers $x$ in $L$), and if $z \parallel_{L \setminus \{x\}} l$, then $l \leq_L w$ for some $w \in L$ such that $w \neq z$ and $w$ lower covers $x$ in $L$, contradicting the join-irreducibility of $x$ in $L$.) Now, if $\inf_L\{y_1,y_2\} \neq x$, then $\inf_L\{y_1,y_2\} \in L \setminus \{x\}$ and $\inf_{L \setminus \{x\}}\{y_1,y_2\} = \inf_L\{y_1,y_2\}$. Also, $y_1,y_2 \in L \setminus \{x\}$ and $x$ being join-irreducible in $L$ together imply that $\sup_L\{y_1,y_2\} \in L \setminus \{x\}$, so we have that $\sup_{L \setminus \{x\}}\{y_1,y_2\} = \sup_L\{y_1, y_2\}$. Thus, a meet and join exist in $L \setminus \{x\}$ for $y_1$ and $y_2$, and $L \setminus \{x\}$ a lattice. If $x$ is meet-irreducible in $L$, then $L\setminus \{x\}$ is a lattice by duality.

\medskip

\noindent Next, we show by contraposition that for all $x \in L$, if the subposet $L \setminus \{x\}$ of $L$ is a lattice, then $x$ is join or meet-irreducible in $L$. If $x=1_L$ and is join-reducible in $L$, or $x=0_L$ and is meet-reducible in $L$, then there are, respectively, two dual atoms from $L$ without a join in $L \setminus \{x\}$, or two atoms from $L$ without a meet in $L \setminus \{x\}$, so $L \setminus \{x\}$ is not a lattice. If $x$ is doubly reducible in $L$, then there exist two elements $y_1$, $y_2$ that upper cover $x$ and two elements $z_1$, $z_2$ that lower cover $x$ in $L$, all distinct. In this case, if $\sup_{L\setminus \{x\}}\{z_1,z_2\}$ exists, then $x = \sup_L\{z_1,z_2\} <_L \sup_{L \setminus \{x\}}\{z_1, z_2\} \leq_L y_1$ and $x = \sup_L\{z_1,z_2\} <_L \sup_{L \setminus \{x\}}\{z_1, z_2\} \leq_L y_2$. Because $y_1 \neq y_2$, we then have that $x <_L \sup_{L \setminus \{x\}}\{z_1, z_2\} <_L y_1$ or $x <_L \sup_{L \setminus \{x\}}\{z_1, z_2\} <_L y_2$, so $y_1$ or $y_2$ does not upper cover $x$ in $L$, a contradiction. Thus, $\sup_{L \setminus \{x\}}\{z_1,z_2\}$ does not exist, and $L \setminus \{x\}$ is not a lattice, as illustrated in the following figure. This completes the proof of Lemma 1.2.

\smallskip

\begin{figure}[H]

\caption*{\textbf{Figure 1.1:} Example corresponding structures in $L$ and (non-lattice) $L \setminus \{x\}$.}

\smallskip

\begin{center}

\begin{tikzpicture}[scale = 0.8]

\node[fill = white] (0) at (-0.25, 3.5) {In $L$};

\node[circle,
    draw = black,
    fill = white] (1) at (-2.75, 2.75) {};

\node[circle,
    draw = black,
    fill = white] (2) at (-1.5, 1.5) {$y_1$};

\node[circle,
    draw = black,
    fill = white] (3) at (1.5, 1.5) {$y_2$};

\node[circle,
    draw = black,
    fill = white] (4) at (-2.75, 0) {};

\node[circle,
    draw = black,
    fill = white] (5) at (0,0) {$x$};

\node[circle,
    draw = black,
    fill = white] (6) at (1.5, 0) {};

\node[circle,
    draw = black,
    fill = white] (7) at (-1.5, -1.5) {$z_1$};

\node[circle,
    draw = black,
    fill = white] (8) at (1.5, -1.5) {$z_2$};

\node[circle,
    draw = black,
    fill = white] (9) at (-2.75, -2.75) {};

\draw (1) -- (2);
\draw (1) -- (4);
\draw (2) -- (5);
\draw (3) -- (5);
\draw (3) -- (6);
\draw (4) -- (9);
\draw (5) -- (7);
\draw (5) -- (8);
\draw (6) -- (8);
\draw (7) -- (9);

\node[fill = white] (0) at (9, 3.5) {In $L \setminus \{x\}$};

\node[circle,
    draw = black,
    fill = white] (1) at (6.5, 2.75) {};

\node[circle,
    draw = black,
    fill = white] (2) at (7.75, 1.5) {$y_1$};

\node[circle,
    draw = black,
    fill = white] (3) at (10.75, 1.5) {$y_2$};

\node[circle,
    draw = black,
    fill = white] (4) at (6.5, 0) {};

\node[circle,
    draw = black,
    fill = white] (5) at (10.75, 0) {};

\node[circle,
    draw = black,
    fill = white] (6) at (7.75, -1.5) {$z_1$};

\node[circle,
    draw = black,
    fill = white] (7) at (10.75, -1.5) {$z_2$};

\node[circle,
    draw = black,
    fill = white] (8) at (6.5, -2.75) {};

\draw (1) -- (4);
\draw (2) -- (1);
\draw (2) -- (6);
\draw (2) -- (7);
\draw (3) -- (5);
\draw (3) -- (6);
\draw (4) -- (8);
\draw (5) -- (7);
\draw (6) -- (8);

\end{tikzpicture}

\end{center}

\end{figure}
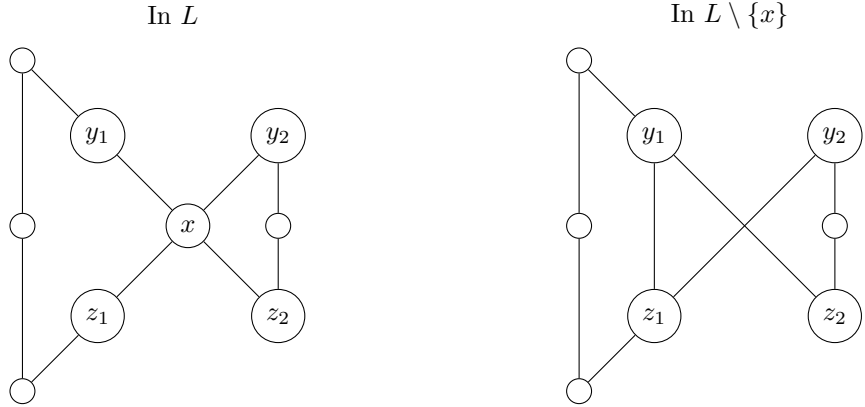

\noindent \textbf{Lemma 1.3.} (\romannumeral 1) \textit{If} $j_1$ \textit{and} $j_2$ \textit{are join-irreducible elements in a finite lattice} $L$\textit{, then} $j_2$ \textit{is join-irreducible in the subposet} $L \setminus \{j_1\}$ \textit{of} $L$; (\romannumeral 2) \textit{If} $m_1$ \textit{and} $m_2$ \textit{are meet-irreducible elements in a finite lattice} $L$\textit{, then} $m_2$ \textit{is meet-irreducible in the subposet} $L \setminus \{m_1\}$ \textit{of} $L$\textit{.}

\medskip

\noindent \textit{Proof.} For part (\romannumeral 1), we denote by $l_1$ and $l_2$ the unique elements that respectively lower cover $j_1$ and $j_2$ in $L$. If $l_2 \neq j_1$, then $l_2$ is the only element that lower covers $j_2$ in $L \setminus \{j_1\}$. Else, $l_2 = j_1$, making $l_1$ the only element that lower covers $l_2$ in $L$. Together with $l_2$ being the only element that lower covers $j_2$ in $L$, this implies that $l_1$ is the only element that lower covers $j_2$ in $L \setminus \{j_1\}$. Thus, $j_2$ is also join-irreducible in $L \setminus \{j_1\}$. Part (\romannumeral 2) follows by duality.

\medskip

We note that the converse of Lemma 1.3 does not hold (for both parts (\romannumeral 1) and (\romannumeral 2)). We consider part (\romannumeral 1). Let $j$ be a join-irreducible element in a finite lattice $L$, and $x$ be an element that upper covers only $j$ and one other element $y$ in $L$ such that $y$ upper covers the unique element $z$ that lower covers $j$ in $L$. Here, the converse of the lemma is not satisfied in that, although $x$ is join-irreducible in the subposet $L \setminus \{j\}$ of $L$, $x$ is not join-irreducible in $L$.

\medskip

\noindent \textbf{Theorem 1.4.} (\romannumeral 1) \textit{If} $J$ \textit{is a set of join-irreducible elements in a finite lattice} $L$\textit{, then the subposet} $L \setminus J$ \textit{of} $L$ \textit{is a lattice}; (\romannumeral 2) \textit{If} $M$ \textit{is a set of meet-irreducible elements in a finite lattice} $L$\textit{, then the subposet} $L \setminus M$ \textit{of} $L$ \textit{is a lattice.}

\medskip

\noindent \textit{Proof.} Again, we consider only part (\romannumeral 1). Without loss of generality, let $J=\{j_1, \cdots, j_{|J|}\}$. If $|J| = 1$, then the subposet $L \setminus J$ of $L$ is a lattice by Lemma 1.2. If $|J|>1$, then we first have that the subposet $L_1 = L \setminus \{j_1\}$ of $L$ is a lattice by Lemma 1.2, and $j_2, \cdots, j_{|J|}$ are join-irreducible in $L_1$ by Lemma 1.3. Then, for all integers $k$ such that $1 \leq k \leq |J|-1$, if the subposet $L_k = L \setminus \{j_1, \cdots, j_k\}$ of $L$ is a lattice and $j_{k+1}, \cdots, j_{|J|}$ are join-irreducible in $L_k$, then the subposet $L_{k+1} = L_k \setminus \{j_{k+1}\}$ of $L_k$ is a lattice by Lemma 1.2 and, when $k<|J|-1$, $j_{k+2}, \cdots, j_{|J|}$ are join-irreducible in $L_{k+1}$ by Lemma 1.3. Thus, we have that the subposet $L \setminus J = L_{|J|} = L_{|J|-1} \setminus \{j_{|J|}\}$ of $L_{|J|-1}$ (and therefore of $L$) is a lattice, completing the proof of Theorem 1.4.

\medskip

Theorem 1.4 does not extend, in general, to a finite lattice $L$ and set $I$ of join and meet-irreducible elements in $L$. A counterexample that demonstrates this is any finite lattice $L$ with $I = \{j, m\}$ such that $j$ is join-irreducible and meet-reducible in $L$, $m$ is meet-irreducible and join-reducible in $L$, and $j$ upper covers $m$ in $L$. In this case, the conclusion of Theorem 1.4 is not satisfied because the subposet $L \setminus I$ of $L$ is not a lattice.

\section*{\large{2. Characterizing a counterexample to Conjecture 1.1 of minimum size}}

We now obtain necessary conditions for any counterexample $\tilde{L}$ to Conjecture 1.1 such that no counterexample $\hat{L}$ exists with $|\hat{L}|<|\tilde{L}|$. We note that $\tilde{L}$ must contain more than two elements.

\medskip

\noindent \textbf{Theorem 2.1.} \textit{Every join-irreducible element} $j$ \textit{in} $\tilde{L}$ \textit{(with the unique element that it upper covers in} $\tilde{L}$ \textit{denoted by} $x$\textit{) lower covers an element} $y$ \textit{in} $\tilde{L}$\textit{, where} $y$ \textit{is not less than or equal to any join-irreducible element and} $y$ \textit{upper covers exactly one other element} $z$ \textit{in} $\tilde{L}$ \textit{such that} $x <_{\tilde{L}} z$\textit{.}

\medskip

\noindent \textit{Proof.} If not, then there exists a join-irreducible element $j$ in $\tilde{L}$ such that any element $y$ that upper covers $j$ in $\tilde{L}$ satisfies at least one of the following criteria:

\medskip

\noindent (\romannumeral 1) $y$ upper covers more than two elements in $\tilde{L}$;

\medskip

\noindent (\romannumeral 2) $y \leq_{\tilde{L}} j^*$ for some join-irreducible $j^*$ in $\tilde{L}$;

\medskip

\noindent (\romannumeral 3) $y$ upper covers $j$ and exactly one other element $z$ in $\tilde{L}$ such that $z$ is incomparable with $x$ in $\tilde{L}$.

\medskip

\noindent We consider the subposet $\hat{L}=\tilde{L} \setminus \{j\}$ of $\tilde{L}$, which is a lattice by Lemma 1.2, and consider any join-irreducible element $\hat{\jmath}$ in $\hat{L}$. First, we assume that $\hat{\jmath}$ is also join-irreducible in $\tilde{L}$. If $\hat{\jmath} <_{\tilde{L}} j$, then $|({\uparrow}\hat{\jmath})_{\hat{L}}| \geq |({\uparrow}j)_{\tilde{L}}| >\frac{|\tilde{L}|}{2}=\frac{|\hat{L}|+1}{2}$. Else, $|({\uparrow}\hat{\jmath})_{\hat{L}}| = |({\uparrow}\hat{\jmath})_{\tilde{L}}| > \frac{|\tilde{L}|}{2}=\frac{|\hat{L}|+1}{2}$. Next, we assume that $\hat{\jmath}$ is join-reducible in $\tilde{L}$. In this case, $\hat{\jmath}$ upper covers $j$ in $\tilde{L}$, and so must satisfy at least one of the three criteria in question. (\romannumeral 1) is not satisfied, as it implies that $\hat{\jmath}$ is join-reducible in $\hat{L}$. (\romannumeral 3) is not satisfied, as it implies that $\hat{\jmath}$ upper covers both $x$ and $z$ in $\hat{L}$, and thus again that $\hat{\jmath}$ is join-reducible in $\hat{L}$. Hence, (\romannumeral 2) is satisfied, and $|({\uparrow}\hat{\jmath})_{\hat{L}}| = |({\uparrow}\hat{\jmath})_{\tilde{L}}| \geq |({\uparrow}j^*)_{\tilde{L}}| > \frac{|\tilde{L}|}{2}=\frac{|\hat{L}|+1}{2}$. Therefore, any join-irreducible element $\hat{\jmath}$ in $\hat{L}$ has $|({\uparrow}\hat{\jmath})_{\hat{L}}| > \frac{|\hat{L}|}{2}$, and $\hat{L}$ is a counterexample to Conjecture 1.1, contradicting the minimality of $|\tilde{L}|$. This completes the proof of Theorem 2.1.

\medskip

In Figure 2.1, $y_1$ satisfies the first criterion from the proof of Theorem 2.1, $y_2$ satisfies the second, and $y_3$ the third. By the theorem, element $j$ of Figure 2.1 must also lower cover an element in $\tilde{L}$ that does not satisfy any of the three criteria.

\medskip

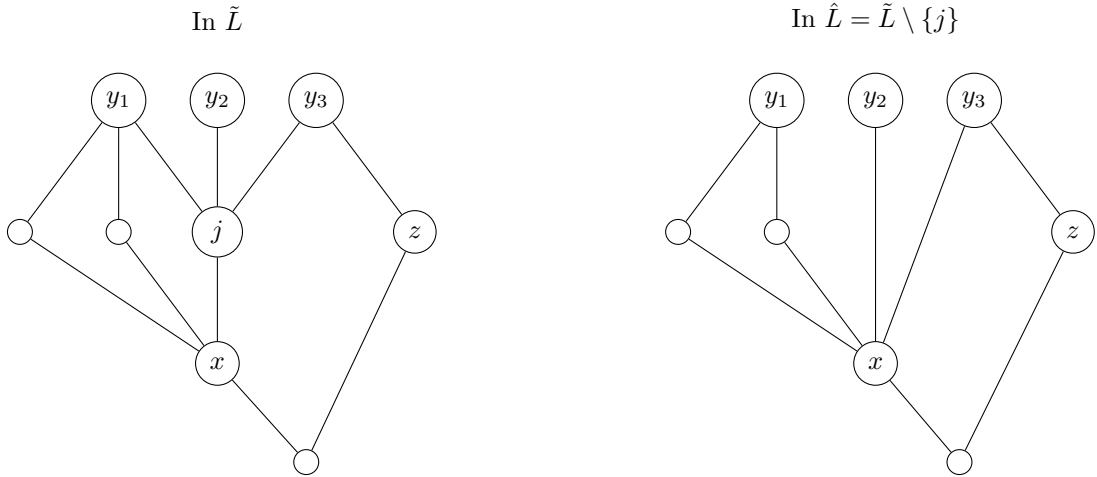
\begin{figure}[H]

\caption*{\textbf{Figure 2.1:} Example corresponding structures in $\tilde{L}$ and $\hat{L}=\tilde{L}\setminus\{j\}$.}

\medskip

\begin{center}

\begin{tikzpicture}[scale = 0.875]

\node[fill = white] (0) at (-1, 6.75) {In $\tilde{L}$};

\node[circle,
    draw = black,
    fill = white] (1) at (-2.5, 5.5) {${y_1}$};

\node[circle,
    draw = black,
    fill = white] (2) at (-1, 5.5) {${y_2}$};

\node[circle,
    draw = black,
    fill = white] (3) at (0.5, 5.5) {${y_3}$};

\node[circle,
    draw = black,
    fill = white] (4) at (-4, 3.5) {};

\node[circle,
    draw = black,
    fill = white] (5) at (-2.5, 3.5) {};

\node[circle,
    draw = black,
    fill = white] (6) at (-1, 3.5) {$j$};

\node[circle,
    draw = black,
    fill = white] (7) at (2, 3.5) {$z$};

\node[circle,
	draw = black,
    fill = white] (8) at (-1, 1.5) {$x$};

\node[circle,
    draw = black,
    fill = white] (9) at (0.35, 0) {};

\draw (1) -- (4);
\draw (1) -- (5);
\draw (1) -- (6); 
\draw (2) -- (6);
\draw (3) -- (6);
\draw (3) -- (7);
\draw (4) -- (8);
\draw (5) -- (8);
\draw (6) -- (8);
\draw (7) -- (9);
\draw (8) -- (9);

\node[fill = white] (0) at (9, 6.75) {In $\hat{L}=\tilde{L}\setminus \{j\}$};

\node[circle,
    draw = black,
    fill = white] (1) at (7.5, 5.5) {$y_1$};

\node[circle,
    draw = black,
    fill = white] (2) at (9, 5.5) {$y_2$};

\node[circle,
    draw = black,
    fill = white] (3) at (10.5, 5.5) {$y_3$};

\node[circle,
    draw = black,
    fill = white] (4) at (6, 3.5) {};

\node[circle,
    draw = black,
    fill = white] (5) at (7.5, 3.5) {};

\node[circle,
    draw = black,
    fill = white] (6) at (12, 3.5) {$z$};

\node[circle,
    draw = black,
    fill = white] (7) at (9, 1.5) {$x$};

\node[circle,
    draw = black,
    fill = white] (8) at (10.275, 0) {};

\draw (1) -- (4);
\draw (1) -- (5);
\draw (2) -- (7);
\draw (3) -- (6);
\draw (3) -- (7);
\draw (4) -- (7);
\draw (5) -- (7);
\draw (6) -- (8);
\draw (7) -- (8);

\end{tikzpicture}

\end{center}

\end{figure}

\noindent \textbf{Corollary 2.2.} $0_{\tilde{L}}$ \textit{is meet-reducible in} $\tilde{L}$\textit{.}

\medskip

\noindent \textit{Proof.} Otherwise, denote by $a$ the unique atom of $\tilde{L}$. For all $x \in \tilde{L} \setminus \{a\}$, $\sup_{\tilde{L}}\{a,x\} \in \{a,x\}$. Because $a$ is join-irreducible in $\tilde{L}$, we have by Theorem 2.1 that there exists an element $b$ upper covering $a$ and exactly one other element $c$ in $\tilde{L}$, making $\sup_{\tilde{L}}\{a,c\} = b \not \in \{a,c\}$, a contradiction.

\medskip

\noindent \textbf{Theorem 2.3.} \textit{No meet-irreducible element is less than a join-irreducible element in} $\tilde{L}$\textit{.}

\medskip

\noindent \textit{Proof.} Otherwise, there exist a meet-irreducible element $m$ and join-irreducible element $j$ in $\tilde{L}$ such that $m <_{\tilde{L}} j$. By Lemma 1.2, the subposet $\hat{L}=\tilde{L} \setminus \{m\}$ of $\tilde{L}$ is a lattice. Consider any join-irreducible element $\hat{\jmath}$ in $\hat{L}$. We first assume that $\hat{\jmath}$ is also join-irreducible in $\tilde{L}$. If $\hat{\jmath} <_{\tilde{L}} m$, then $\hat{\jmath} <_{\tilde{L}} j$ and $|({\uparrow}\hat{\jmath})_{\hat{L}}|=|({\uparrow}\hat{\jmath})_{\tilde{L}}|-1 > |({\uparrow}j)_{\tilde{L}}| > \frac{|\tilde{L}|}{2} = \frac{|\hat{L}|+1}{2}$. Else, $|({\uparrow}\hat{\jmath})_{\hat{L}}|=|({\uparrow}\hat{\jmath})_{\tilde{L}}| > \frac{|\tilde{L}|}{2}=\frac{|\hat{L}|+1}{2}$. Next, we instead assume that $\hat{\jmath}$ is join-reducible in $\tilde{L}$, implying that $\hat{\jmath}$ upper covers $m$ in $\tilde{L}$. In this case, $\hat{\jmath} <_{\tilde{L}} j$, and we have that $|({\uparrow}\hat{\jmath})_{\hat{L}}|=|({\uparrow}\hat{\jmath})_{\tilde{L}}| > |({\uparrow}j)_{\tilde{L}}| > \frac{|\tilde{L}|}{2}=\frac{|\hat{L}|+1}{2}$. In both cases, $|({\uparrow}\hat{\jmath})_{\hat{L}}| > \frac{|\hat{L}|}{2}$, and so $\hat{L}$ is a counterexample to Conjecture 1.1 that is smaller than $\tilde{L}$. This contradicts the minimality of $|\tilde{L}|$, completing the proof of Theorem 2.3.

\medskip

\noindent \textbf{Corollary 2.4.} $1_{\tilde{L}}$ \textit{is join-reducible in} $\tilde{L}$\textit{.}

\medskip

\noindent \textit{Proof.} If not, then $1_{\tilde{L}}$ is join-irreducible in $\tilde{L}$, upper covering exactly one dual atom. This contradicts Theorem 2.3, as all dual atoms of a lattice are meet-irreducible.

\medskip

\noindent \textbf{Lemma 2.5.} \textit{If} $x$ \textit{is doubly irreducible in} $\tilde{L}$\textit{, then} $|({\uparrow}x)_{\tilde{L}}| = \frac{|\tilde{L}|+1}{2}$\textit{.}

\medskip

\noindent \textit{Proof.} Assume otherwise, i.e. that there is a doubly irreducible element $x$ in $\tilde{L}$ such that $|({\uparrow}x)_{\tilde{L}}| > \frac{|\tilde{L}|+1}{2}$. By Lemma 1.2, the subposet $\hat{L}=\tilde{L} \setminus \{x\}$ of $\tilde{L}$ is a lattice. For any join-irreducible element $\hat{\jmath}$ in $\hat{L}$, we first assume that $\hat{\jmath}$ is also join-irreducible in $\tilde{L}$. If $\hat{\jmath} <_{\tilde{L}} x$, then $|({\uparrow}\hat{\jmath})_{\hat{L}}|=|({\uparrow}\hat{\jmath})_{\tilde{L}}|-1 \geq |({\uparrow}x)_{\tilde{L}}| > \frac{|\tilde{L}|+1}{2} = \frac{|\hat{L}|+2}{2}$. Else, $|({\uparrow}\hat{\jmath})_{\hat{L}}|=|({\uparrow}\hat{\jmath})_{\tilde{L}}| > \frac{|\tilde{L}|}{2} = \frac{|\hat{L}|+1}{2}$. Now we assume that $\hat{\jmath}$ is join-reducible in $\tilde{L}$. In this case, $\hat{\jmath}$ upper covers $x$ in $\tilde{L}$, and $|({\uparrow}\hat{\jmath})_{\hat{L}}|=|({\uparrow}\hat{\jmath})_{\tilde{L}}| = |({\uparrow}x)_{\tilde{L}}|-1 > \frac{|\tilde{L}|+1}{2}-1 = \frac{|\hat{L}|}{2}$. In both cases, $|({\uparrow}\hat{\jmath})_{\hat{L}}| > \frac{|\hat{L}|}{2}$. Therefore, $\hat{L}$ is a counterexample to Conjecture 1.1, contradicting the minimality of $|\tilde{L}|$. This completes the proof of Lemma 2.5.

\medskip

\noindent \textbf{Theorem 2.6.} \textit{There is at most one doubly irreducible element in }$\tilde{L}$\textit{.}

\medskip

\noindent \textit{Proof.} If not, then there exist two distinct elements $x$ and $y$ that are doubly irreducible in $\tilde{L}$. We first assume that $x$ and $y$ are comparable in $\tilde{L}$, where $x <_{\tilde{L}} y$ without loss of generality. Then $|({\uparrow}x)_{\tilde{L}}|>|({\uparrow}y)_{\tilde{L}}|$, contradicting Lemma 2.5 which requires that $|({\uparrow}x)_{\tilde{L}}|=|({\uparrow}y)_{\tilde{L}}|=\frac{|\tilde{L}|+1}{2}$. Next, we assume that $x$ and $y$ are incomparable in $\tilde{L}$. We consider the subposet $\hat{L}=L \setminus \{x, y\}$ of $\tilde{L}$, which is itself a lattice by Theorem 1.4. Now, consider any join-irreducible element $\hat{\jmath}$ in $\hat{L}$. First, we assume that $\hat{\jmath}$ is also join-irreducible in $\tilde{L}$. If either $\hat{\jmath} <_{\tilde{L}} x$ or $\hat{\jmath} <_{\tilde{L}} y$, then $|({\uparrow}\hat{\jmath})_{\hat{L}}|=|({\uparrow}\hat{\jmath})_{\tilde{L}}|-1 > \frac{|\tilde{L}|+1}{2}-1=\frac{|\hat{L}|+1}{2}$. If neither $\hat{\jmath} <_{\tilde{L}} x$ nor $\hat{\jmath} <_{\tilde{L}} y$, then $|({\uparrow}\hat{\jmath})_{\hat{L}}|=|({\uparrow}\hat{\jmath})_{\tilde{L}}| > \frac{|\tilde{L}|}{2}=\frac{|\hat{L}|+2}{2}$. If both $\hat{\jmath} <_{\tilde{L}} x$ and $\hat{\jmath} <_{\tilde{L}} y$, then $|({\uparrow}\hat{\jmath})_{\hat{L}}|=|({\uparrow}\hat{\jmath})_{\tilde{L}}|-2 > \frac{|\tilde{L}|+3}{2}-2=\frac{|\hat{L}|+1}{2}$. Next, we assume that $\hat{\jmath}$ is join-reducible in $\tilde{L}$. In this case, at least one of $x$ or $y$ lower covers $\hat{\jmath}$ in $\tilde{L}$, and $|({\uparrow}\hat{\jmath})_{\hat{L}}|=|({\uparrow}\hat{\jmath})_{\tilde{L}}| = \frac{|\tilde{L}|+1}{2}-1 = \frac{|\hat{L}|+1}{2}$. Therefore, in both cases $|({\uparrow}\hat{\jmath})_{\hat{L}}| > \frac{|\hat{L}|}{2}$, making $\hat{L}$ a counterexample to Conjecture 1.1. This contradicts the minimality of $|\tilde{L}|$, completing the proof of Theorem 2.6.

\medskip

In [4], two incomparable doubly irreducible elements were also removed from a lattice, in that case with respect to the class of dismantlable lattices, in order to inductively prove that Conjecture 1.1 holds for all lattices therein. The \textit{length} $\ell(L)$ of a lattice $L$ is one less than the maximum size of a chain in $L$. Theorem 2.6 can be paired with a result of Rival (see Theorem 1 of [9]) to prove the following:

\medskip

\noindent \textbf{Theorem 2.7.} \textit{If} $j$ \textit{is a join-irreducible element in} $\tilde{L}$\textit{, then} $|({\uparrow}j)_{\tilde{L}}| > \ell(\tilde{L})$\textit{.}

\medskip

\noindent \textit{Proof.} In [9], it was shown that $|L| \geq 2(\ell(L)+1)-|\texttt{Irr}(L)|$ for any finite lattice $L$, where $\texttt{Irr}(L)$ is the set of doubly irreducible elements in $L$. We note that $0_L$ and $1_L$ are considered by definition in [9] to be, respectively, join and meet-irreducible in a finite lattice $L$, which is not the case in the present work. However, the set of elements that are doubly irreducible in $\tilde{L}$ is not affected by this difference because $0_{\tilde{L}}$ is meet-reducible in $\tilde{L}$ by Corollary 2.2, and $1_{\tilde{L}}$ is join-reducible in $\tilde{L}$ by Corollary 2.4, implying that, whether or not $0_{\tilde{L}}$ and $1_{\tilde{L}}$ are respectively considered join and meet-irreducible in $\tilde{L}$, neither are considered doubly irreducible. Now, by Theorem 2.6, $\texttt{Irr}(\tilde{L}) \leq 1$. It follows that $|\tilde{L}| \geq 2(\ell(\tilde{L})+1)-1$, and so $\frac{|\tilde{L}|}{2} > \ell(\tilde{L})$. To complete the proof of Theorem 2.7, we recall that any join-irreducible element $j$ in $\tilde{L}$ has $|({\uparrow}j)_{\tilde{L}}| > \frac{|\tilde{L}|}{2}$, and thus also $|({\uparrow}j)_{\tilde{L}}| > \ell(\tilde{L})$.

\medskip

\noindent \textbf{Corollary 2.8.} \textit{Any doubly irreducible element is less than at least one doubly reducible element in} $\tilde{L}$\textit{.}

\medskip

\noindent \textit{Proof.} By Theorem 2.3, every element greater than a doubly irreducible element $x$ in $\tilde{L}$ is join-reducible in $\tilde{L}$, so we need only to show that one such element is also meet-reducible in $\tilde{L}$. We assume that every element $y \in ({\uparrow}x)_{\tilde{L}} \setminus \{1_{\tilde{L}}\}$ is meet-irreducible in $\tilde{L}$. Then there is a unique maximal chain $C$ in $\tilde{L}$ such that $0_C=x$ and $1_C=1_{\tilde{{L}}}$, and an element of $\tilde{L}$ is greater than $x$ in $\tilde{L}$ if and only if it is greater than $x$ in $C$. Further, $0_{\tilde{L}} <_{\tilde{L}} x$ implies that $C$ is a subposet of a chain $C_0$ in $\tilde{L}$ such that $0_{\tilde{L}} \in C_0$. It follows that $|({\uparrow}x)_{\tilde{L}}| = |C| < |C_0| \leq \ell(\tilde{L})+1$. Therefore, $|({\uparrow}x)_{\tilde{L}}| \leq \ell(\tilde{L})$, contradicting Theorem 2.7.

\medskip

We observe that every element $x \in \tilde{L} \setminus \{0_{\tilde{L}}, 1_{\tilde{L}}\}$ is comparable with at least four elements, excluding itself, in $\tilde{L}$. These include $0_{\tilde{L}}$ and $1_{\tilde{L}}$, as well as two other distinct elements from $\tilde{L} \setminus \{x\}$. (Otherwise, $x$ is doubly irreducible in $\tilde{L}$, and either lower covers or is equal to a dual atom in $\tilde{L}$, contradicting Corollary 2.8.)

\medskip

\noindent \textbf{Theorem 2.9.} \textit{If} $M$ \textit{is a nonempty set of elements that are meet-irreducible in} $\tilde{L}$\textit{, then there exists a join-irreducible element} $j$ \textit{in} $\tilde{L}$ \textit{such that} $|({\uparrow}j)_{\tilde{L}} \cap M| > \frac{|M|}{2}$\textit{.}

\medskip

\noindent \textit{Proof.} Otherwise, there exists a nonempty set $M$ of meet-irreducible elements in $\tilde{L}$ such that $|({\uparrow}j)_{\tilde{L}} \cap M| \leq \frac{|M|}{2}$ for every join-irreducible element $j$ in $\tilde{L}$. We consider the subposet $\hat{L}=\tilde{L} \setminus M$ of $\tilde{L}$, which is a lattice by Theorem 1.4, and consider any join-irreducible element $\hat{\jmath}$ in $\hat{L}$. If $\hat{\jmath}$ is also join-irreducible in $\tilde{L}$, then $|({\uparrow}\hat{\jmath})_{\hat{L}}| > \frac{|\tilde{L}|}{2}-|({\uparrow}\hat{\jmath})_{\tilde{L}} \cap M|=\frac{|\hat{L}|+|M|}{2}-|({\uparrow}\hat{\jmath})_{\tilde{L}} \cap M|=\frac{|\hat{L}|}{2}+(\frac{|M|}{2}-|({\uparrow}\hat{\jmath})_{\tilde{L}} \cap M|) \geq \frac{|\hat{L}|}{2}$. If $\hat{\jmath}$ is join-reducible in $\tilde{L}$, then there exists $m_1 \in M$ that lower covers $\hat{\jmath}$ in $\tilde{L}$. It follows that $\{m_1\} \subseteq C=\{m_1, \cdots , m_{|C|}\} \subseteq M$, where $C$ is any chain of maximum size in $\tilde{L}$ such that $1 \leq i < |C|$ implies that $m_i$ upper covers $m_{i+1}$ in $\tilde{L}$. We first assume that $m_{|C|}$ is join-irreducible in $\tilde{L}$, implying that $|({\uparrow}m_{|C|})_{\tilde{L}}| > \frac{|\tilde{L}|}{2}$. Then, based on the definition of $C$, we have that $|({\uparrow}\hat{\jmath})_{\tilde{L}}| = |({\uparrow}m_{|C|})_{\tilde{L}}|-|C| > \frac{|\tilde{L}|-2|C|}{2}$, and so $|({\uparrow}\hat{\jmath})_{\hat{L}}| = |({\uparrow}\hat{\jmath})_{\tilde{L}}|-|({\uparrow}\hat{\jmath})_{\tilde{L}} \cap M| > \frac{|\hat{L}|+|M|-2(|C|+|({\uparrow}\hat{\jmath})_{\tilde{L}} \cap M|)}{2}$. If $|C|+|({\uparrow}\hat{\jmath})_{\tilde{L}} \cap M| > \frac{|M|}{2}$, then $|({\uparrow}m_{|C|})_{\tilde{L}} \cap M| > \frac{|M|}{2}$, as $C \cup (({\uparrow}\hat{\jmath})_{\tilde{L}} \cap M) = ({\uparrow}m_{|C|})_{\tilde{L}} \cap M$ and $C \cap (({\uparrow}\hat{\jmath})_{\tilde{L}} \cap M) = \emptyset$. However, because $m_{|C|}$ is join-irreducible in $\tilde{L}$, we also have that $|({\uparrow}m_{|C|})_{\tilde{L}} \cap M| \leq \frac{|M|}{2}$, a contradiction. Thus, $|C|+|({\uparrow}\hat{\jmath})_{\tilde{L}} \cap M| \leq \frac{|M|}{2}$, making $|({\uparrow}\hat{\jmath})_{\hat{L}}| > \frac{|\hat{L}|}{2}$. Now, we assume that $m_{|C|}$ is join-reducible in $\tilde{L}$, upper covering distinct elements $x,y \not \in M$. In this case, because $\hat{\jmath}$ is join-irreducible in $\hat{L}$, we have that $|C|>1$, and there exists an element $z \not \in C$ less than $\hat{\jmath}$ in $\tilde{L}$ that upper covers $m_k$ in $\tilde{L}$ for some $k \in \{2, \cdots, |C|\}$. This contradicts the meet-irreducibility of $m_k$ in $\tilde{L}$. In conclusion, whether or not $\hat{\jmath}$ is join-irreducible in $\tilde{L}$, we have that $|({\uparrow}\hat{\jmath})_{\hat{L}}| > \frac{|\hat{L}|}{2}$. Then $\hat{L}$ is a counterexample to Conjecture 1.1, contradicting the minimality of $|\tilde{L}|$. This completes the proof of Theorem 2.9.

\medskip

\noindent \textbf{Corollary 2.10.} \textit{There is a join-irreducible element} $j$ \textit{in} $\tilde{L}$ \textit{that is less than or equal to more than half of the meet-irreducible elements in} $\tilde{L}$.

\medskip

\noindent \textit{Proof.} We set $M$ of Theorem 2.9 equal to the set of all meet-irreducible elements in $\tilde{L}$.

\medskip

\noindent \textbf{Corollary 2.11.} \textit{For any two meet-irreducible elements} $m_1$ \textit{and} $m_2$ \textit{in} $\tilde{L}$\textit{, there exists a join-irreducible element} $j$ \textit{in} $\tilde{L}$ \textit{such that} $j \leq_{\tilde{L}} m_1$ \textit{and} $j \leq_{\tilde{L}} m_2$\textit{.}

\medskip

\noindent \textit{Proof.} We set $M$ of Theorem 2.9 equal to $\{m_1,m_2\}$.

\medskip

By Corollary 2.11, any meet-irreducible atom is less than all other meet-irreducible elements in $\tilde{L}$.

\medskip

\noindent \textbf{Theorem 2.12.} \textit{For every meet-irreducible element} $m$ \textit{in} $\tilde{L}$\textit{, there exists a join-irreducible element} $j$ \textit{in} $\tilde{L}$ \textit{such that} $j \leq_{\tilde{L}} m$ \textit{and} $|({\uparrow}j)_{\tilde{L}}|=\frac{|\tilde{L}|+1}{2}$.

\medskip

\noindent \textit{Proof.} Otherwise, there exists a meet-irreducible element $m$ in $\tilde{L}$ such that any join-irreducible element $j$ in $\tilde{L}$ with $j \leq_{\tilde{L}} m$ has $|({\uparrow}j)_{\tilde{L}}| > \frac{|\tilde{L}|+1}{2}$. We consider the subposet $\hat{L}=\tilde{L} \setminus \{m\}$ of $\tilde{L}$, which is a lattice by Lemma 1.2. Any join-irreducible element $\hat{\jmath}$ in $\hat{L}$ is also join-irreducible in $\tilde{L}$. (If not, then $\hat{\jmath}$ upper covers $m$ in $\tilde{L}$, and because $\hat{\jmath}$ is join-irreducible in $\hat{L}$, we have that $m$ is join-irreducible in $\tilde{L}$. It follows that $|({\uparrow}m)_{\tilde{L}}| > \frac{|\tilde{L}|+1}{2}$, yet $|({\uparrow}m)_{\tilde{L}}| = \frac{|\tilde{L}|+1}{2}$ by Lemma 2.5, a contradiction.) If $\hat{\jmath} <_{\tilde{L}} m$, then $|({\uparrow}\hat{\jmath})_{\hat{L}}| = |({\uparrow}\hat{\jmath})_{\tilde{L}}| - 1 > \frac{|\tilde{L}|+1}{2}-1 = \frac{|\hat{L}|}{2}$. Else, $|({\uparrow}\hat{\jmath})_{\hat{L}}| = |({\uparrow}\hat{\jmath})_{\tilde{L}}| > \frac{|\tilde{L}|}{2} = \frac{|\hat{L}|+1}{2}$. Therefore, $\hat{L}$ is a counterexample to Conjecture 1.1. This contradicts the minimality of $|\tilde{L}|$, completing the proof of Theorem 2.12.

\medskip

\noindent \textbf{Theorem 2.13.} \textit{For all} $x \in \tilde{L} \setminus \{0_{\tilde{L}}, 1_{\tilde{L}}\}$\textit{, the subposet} $({\parallel}x)_{\tilde{L}}$ \textit{of} $\tilde{L}$ \textit{is not a chain in} $\tilde{L}$\textit{.}

\medskip

\noindent \textit{Proof.} Assume otherwise. Then there exists an element $x \in \tilde{L} \setminus \{0_{\tilde{L}}, 1_{\tilde{L}}\}$ such that $({\parallel}x)_{\tilde{L}}=\emptyset$ or $({\parallel}x)_{\tilde{L}}=\{c_1, \cdots, c_n \}$, where $n=|({\parallel}x)_{\tilde{L}}|$ and $1 \leq i \leq j \leq n$ implies that $c_j \leq_{\tilde{L}} c_i$ without loss of generality. No element from $({\parallel}x)_{\tilde{L}}$ can upper cover an element from $({\uparrow}x)_{\tilde{L}}$ or lower cover an element from $({\downarrow}x)_{\tilde{L}}$ in $\tilde{L}$. Thus, if $|({\parallel}x)_{\tilde{L}}|>1$ and $1 \leq i < n$, then $c_i$ upper covers $c_{i+1}$ in $\tilde{L}$. Further, if $|({\parallel}x)_{\tilde{L}}|>1$ and $1 \leq i < j \leq n$, and two elements $u_i$ and $u_j$ ($l_i$ and $l_j$) from $({\uparrow}x)_{\tilde{L}}$ (from $({\downarrow}x)_{\tilde{L}}$) respectively upper cover (lower cover) $c_i$ and $c_j$ in $\tilde{L}$, then $u_j <_{\tilde{L}} u_i$ ($l_j <_{\tilde{L}} l_i$). We consider the subposet $\hat{L} = ({\uparrow}x)_{\tilde{L}}$ of $\tilde{L}$. For all $y_1, y_2 \in \hat{L}$, $\sup_{\tilde{L}}\{y_1, y_2\} \in \hat{L}$ because $x \leq_{\tilde{L}} y_1 \leq_{\tilde{L}} \sup_{\tilde{L}}\{y_1,y_2\}$. Also, $\inf_{\tilde{L}}\{y_1, y_2\} \in \hat{L}$ because $x \leq_{\tilde{L}} y_1$ and $x \leq_{\tilde{L}} y_2$ together imply that $x \leq_{\tilde{L}} \inf_{\tilde{L}}\{y_1, y_2\}$. Thus, $\hat{L}$ is a lattice with $\sup_{\hat{L}}\{y_1, y_2\} = \sup_{\tilde{L}}\{y_1,y_2\}$ and $\inf_{\hat{L}}\{y_1, y_2\} = \inf_{\tilde{L}}\{y_1,y_2\}$. Now, consider any join-irreducible element $\hat{\jmath}$ in $\hat{L}$. If $\hat{\jmath}$ is also join-irreducible in $\tilde{L}$, then $|({\uparrow}\hat{\jmath})_{\hat{L}}|=|({\uparrow}\hat{\jmath})_{\tilde{L}}| > \frac{|\tilde{L}|}{2} = \frac{|\hat{L}|+|\tilde{L} \setminus \hat{L}|}{2} > \frac{|\hat{L}|}{2}$. If $\hat{\jmath}$ is join-reducible in $\tilde{L}$, then $\hat{\jmath}$ upper covers an element from $({\parallel}x)_{\tilde{L}}$ in $\tilde{L}$, and we denote the unique element from $({\uparrow}x)_{\tilde{L}}$ that upper covers $c_1$ in $\tilde{L}$ by $u_1$. We first assume that $n=1$, or that $n > 1$ and $c_1$ only upper covers $c_2$ in $\tilde{L}$. In either case, $c_1$ is doubly irreducible in $\tilde{L}$, with $|({\uparrow}c_1)_{\tilde{{L}}}|=\frac{|\tilde{L}|+1}{2}$ by Lemma 2.5, and $\hat{\jmath} \leq_{\tilde{L}} u_1$ implies that $|({\uparrow}\hat{\jmath})_{\tilde{{L}}}| \geq |({\uparrow}u_1)_{\tilde{{L}}}| = |({\uparrow}c_1)_{\tilde{{L}}}|-1=\frac{|\tilde{L}|-1}{2} = \frac{|\hat{L}| + |\tilde{L} \setminus \hat{L}| - 1}{2}$. It follows that $|({\uparrow}\hat{\jmath})_{\tilde{{L}}}| > \frac{|\hat{L}|}{2}$, as $\{0_{\tilde{L}}, c_1\} \subseteq \tilde{L} \setminus \hat{L}$ implies that $|\tilde{L} \setminus \hat{L}| \geq 2$. Next, we assume that $n > 1$ and $c_1$ upper covers an element $l_1 \in ({\downarrow}x)_{\tilde{L}}$ in $\tilde{L}$. This implies that every element of $({\parallel}x)_{\tilde{L}}$ is meet-irreducible in $\tilde{L}$. (Otherwise, there exists $m \in \{2, \dots, n\}$ such that $c_m$ lower covers some $u_m \in ({\uparrow}x)_{\tilde{L}}$ in $\tilde{L}$, and we have that $c_m <_{\tilde{L}} c_1$ and $u_m$ upper covers $c_m$ in $\tilde{L}$, yet $l_1 <_{\tilde{L}} u_m$ and $c_1$ upper covers $l_1$ in $\tilde{L}$, contradicting the existence of $\inf_{\tilde{L}}\{c_1, u_m\}$.) It follows that $\hat{\jmath}=u_1$, and $c_n$ is doubly irreducible in $\tilde{L}$. By Theorem 2.6, we then have that all elements $c_1, \cdots, c_{n-1}$ are join-reducible in $\tilde{L}$. Thus, for every $i \in \{1,\cdots, n\}$, $c_i$ upper covers an element $l_i \in ({\downarrow}x)_{\tilde{L}}$ in $\tilde{L}$, and $1 \leq j < k \leq n$ implies that $l_k <_{\tilde{L}} l_j$, as illustrated in Figure 2.2. Then $|({\uparrow}\hat{\jmath})_{\hat{L}}| = |({\uparrow}\hat{\jmath})_{\tilde{L}}| = |({\uparrow}c_n)_{\tilde{L}}|-n = \frac{|\tilde{L}|+1-2n}{2}=\frac{|\hat{L}|+|\tilde{L} \setminus \hat{L}|+1-2n}{2}$, making $|({\uparrow}\hat{\jmath})_{\hat{L}}| > \frac{|\hat{L}|}{2}$ because $\{c_1, l_1, \cdots, c_n, l_n \} \subseteq \tilde{L} \setminus \hat{L}$ implies that $|\tilde{L} \setminus \hat{L}| \geq 2n$. Therefore, $|({\uparrow}\hat{\jmath})_{\hat{L}}| > \frac{|\hat{L}|}{2}$ whether or not $\hat{\jmath}$ is join-irreducible in $\tilde{L}$, and $\hat{L}$ is a counterexample to Conjecture 1.1. This contradicts the minimality of $|\tilde{L}|$, completing the proof of Theorem 2.13.

\medskip

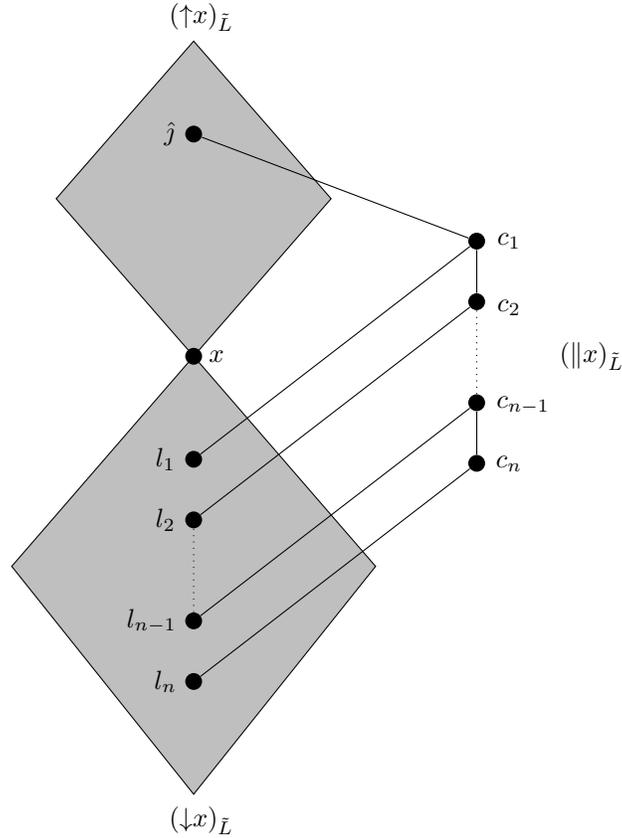
\begin{figure}[H]

\caption*{\textbf{Figure 2.2:} $\hat{L}=({\uparrow}x)_{\tilde{L}}$ must also be a counterexample to Conjecture 1.1.}

\medskip

\begin{center}

\begin{tikzpicture}[scale = 1.075]

\node[fill = white] at (0.1,-2.05) {$({\uparrow}x)_{\tilde{L}}$};
\node[fill = white] at (0.1,-12) {$({\downarrow}x)_{\tilde{L}}$};
\node[fill = white] at (4.925,-6.25) {$({\parallel}x)_{\tilde{L}}$};

\draw [fill = lightgray] (0,-2.35) -- (-1.7,-4.3) -- (0,-6.25) -- (1.7,-4.3) -- cycle;

\draw [fill = lightgray] (0,-6.25) -- (-2.25,-8.85) -- (0,-11.67) -- (2.25,-8.85) -- cycle;

\node at (-0.3,-3.5) {$\hat{\jmath}$};
\node[scale = 0.675, circle, fill = black] (0) at (0,-3.5) {};

\node at (3.9,-4.825) {$c_1$};
\node[scale = 0.675, circle, fill = black] (1a) at (3.5,-4.825) {};

\node at (3.9,-5.65) {$c_2$};
\node[scale = 0.675, circle, fill = black] (1b) at (3.5,-5.575) {};

\node at (4.075,-6.855) {$c_{n-1}$};
\node[scale = 0.675, circle, fill = black] (1c) at (3.5,-6.825) {};

\node at (3.9,-7.575) {$c_n$};
\node[scale = 0.675, circle, fill = black] (1d) at (3.5,-7.575) {};

\draw (0)--(1a);

\node at (0.285,-6.25) {$x$};
\node[scale = 0.675, circle, fill = black] at (0,-6.25) {};

\draw (1a)--(1b);
\draw[dotted] (1b)--(1c);
\draw (1c)--(1d);

\node at (-0.35,-7.525) {$l_1$};
\node[scale = 0.675, circle, fill = black] (2a) at (0,-7.525) {};

\node at (-0.35,-8.275) {$l_2$};
\node[scale = 0.675, circle, fill = black] (2b) at (0,-8.275) {};

\node at (-0.52,-9.525) {$l_{n-1}$};
\node[scale = 0.675, circle, fill = black] (2c) at (0,-9.525) {};

\node at (-0.35,-10.275) {$l_n$};
\node[scale = 0.675, circle, fill = black] (2d) at (0,-10.275) {};

\draw[dotted] (2b)--(2c);

\draw (1a)--(2a);   
\draw (1b)--(2b);   
\draw (1c)--(2c);
\draw (1d)--(2d);

\end{tikzpicture}

\end{center}

\end{figure}

\noindent \textbf{Corollary 2.14.} \textit{Every element} $x \in \tilde{L} \setminus \{0_{\tilde{L}}, 1_{\tilde{L}}\}$ \textit{is incomparable with at least three elements in} $\tilde{L}$\textit{.}

\medskip

\noindent \textit{Proof.} Otherwise, there exists an element $x \in \tilde{L} \setminus \{0_{\tilde{L}}, 1_{\tilde{L}}\}$ such that $|({\parallel}x)_{\tilde{L}}| \in \{0,1,2\}$. If $|({\parallel}x)_{\tilde{L}}| \in \{0,1\}$, then $({\parallel}x)_{\tilde{L}}$ is a chain in $\tilde{L}$, contradicting Theorem 2.13. If $|({\parallel}x)_{\tilde{L}}| = 2$, then $({\parallel}x)_{\tilde{L}}$ is either a chain in $\tilde{L}$, contradicting Theorem 2.13, or both elements of $({\parallel}x)_{\tilde{L}}$ are doubly irreducible in $\tilde{L}$, contradicting Theorem 2.6. Therefore, $|({\parallel}x)_{\tilde{L}}| \geq 3$.

\medskip

By similar reasoning, greater lower bounds can be proved for $|({\parallel}x)_{\tilde{L}}|$. Further, the argument of Theorem 2.13 can be extended to show that the subposet $({\parallel}x)_{\tilde{L}}$ of $\tilde{L}$ is not of the form $C_{(k)}=\bigcup_{1 \leq i \leq k}C_i$, where $\{C_i\}_{1 \leq j \leq k}$ is a sequence of $k$ disjoint chains in $\tilde{L}$ such that $x \leq_{C_{(k)}} y$ implies that $x,y \in C_j$ for some $j \in \{1, \cdots, k\}$.

\medskip

\noindent \textbf{Theorem 2.15.} \textit{If} $L'$ \textit{is a subposet of} $\tilde{L}$ \textit{and itself a lattice satisfying at least one of the following criteria}:

\smallskip

\noindent (\romannumeral 1) $3 < |L'| < 8$;

\smallskip

\noindent (\romannumeral 2) $2 < |L'| < |\tilde{L}|-2$ \textit{and} $d \leq_{\tilde{L}} 1_{L'}$ \textit{for some dual atom} $d$ \textit{in} $\tilde{L}$:

\smallskip

\noindent \textit{Then there exists an element} $x \in L' \setminus \{0_{L'}, 1_{L'}\}$ \textit{that upper or lower covers some element} $y \in \tilde{L} \setminus L'$ \textit{in} $\tilde{L}$\textit{.}

\medskip

\noindent \textit{Proof.} Otherwise, there exists a subposet $L'$ of $\tilde{L}$ that is itself a lattice satisfying (\romannumeral 1) or (\romannumeral 2) such that no element $x \in L' \setminus \{0_{L'}, 1_{L'}\}$ upper or lower covers any element $y \in \tilde{L} \setminus L'$ in $\tilde{L}$. First, we assume that $L'$ satisfies (\romannumeral 1), i.e. that $3 < |L'| < 8$. Provided by Kyuno in [5] are $2$, $5$, $15$, and $53$ Hasse diagrams for all unlabeled lattices $L$ with, respectively, four, five, six, and seven elements, matching the numbers of unlabeled lattices obtained also in [3]. We observe from these diagrams that all lattices $L$ with $3 < |L| < 8$ have at least two doubly irreducible elements. Therefore, there exist two doubly irreducible elements $x_1$ and $x_2$ in $L'$. It follows that $x_1$ and $x_2$ are also doubly irreducible in $\tilde{L}$, as neither $x_1$ nor $x_2$ upper or lower covers any element $x_3 \in \tilde{L} \setminus L'$ in $\tilde{L}$. This contradicts Theorem 2.6, and so $L'$ does not satisfy (\romannumeral 1).

\medskip

\noindent Next, we assume that $L'$ satisfies (\romannumeral 2), i.e. that $2 < |L'| < |\tilde{L}|-2$ and $d \leq_{\tilde{L}} 1_{L'}$ for some dual atom $d$ in $\tilde{L}$. If $1_{L'}$ is join-reducible in $L'$, then we let $\hat{L}$ be the lattice $L'$ and consider any join-irreducible element $\hat{\jmath}$ in $\hat{L}$. For all $x \in \hat{L} \setminus \{0_{\hat{L}}, 1_{\hat{L}}\}$, if $x$ upper covers an element $y$ in $\tilde{L}$, then $y$ belongs to $\hat{L}$. This implies that $\hat{\jmath}$ is also join-irreducible in $\tilde{L}$. Also, if $x$ lower covers an element $y$ in $\tilde{L}$, then $y$ belongs to $\hat{L}$. It follows that $|({\uparrow}\hat{\jmath})_{\hat{L}}| \geq |({\uparrow}\hat{\jmath})_{\tilde{L}}|-1$, as $d \leq_{\tilde{L}} 1_{L'}$. Hence, $|({\uparrow}\hat{\jmath})_{\hat{L}}| > \frac{|\tilde{L}|}{2}-1=\frac{|\hat{L}|+|\tilde{L} \setminus \hat{L}|-2}{2}$, and we have that $|({\uparrow}\hat{\jmath})_{\hat{L}}| > \frac{|\hat{L}|}{2}$, as $|L'| < |\tilde{L}|-2$ and $\hat{L}=L'$ together imply that $|\tilde{L} \setminus \hat{L}| \geq 3$. Now, if $1_{L'}$ is join-irreducible in $L'$, then we let $\hat{L}$ be the subposet $L' \setminus \{1_{L'}\}$ of $L'$, which is a lattice by Lemma 1.2, and again consider any join-irreducible element $\hat{\jmath}$ in $\hat{L}$. In this case, $\hat{\jmath}$ is join-irreducible in $\tilde{L}$ and $|({\uparrow}\hat{\jmath})_{\hat{L}}| \geq |({\uparrow}\hat{\jmath})_{\tilde{L}}|-2$. We have that $|({\uparrow}\hat{\jmath})_{\hat{L}}| > \frac{|\tilde{L}|}{2}-2=\frac{|\hat{L}|+|\tilde{L} \setminus \hat{L}|-4}{2}$, and $|({\uparrow}\hat{\jmath})_{\hat{L}}| > \frac{|\hat{L}|}{2}$, as $|L'| < |\tilde{L}|-2$ and $\hat{L}=L' \setminus \{1_{L'}\}$ imply that $|\tilde{L} \setminus \hat{L}| \geq 4$. Therefore, in both cases $\hat{L}$ is a counterexample to Conjecture 1.1. This contradicts the minimality of $|\tilde{L}|$, and so $L'$ does not satsify (\romannumeral 2). Thus, $L'$ satisfies neither (\romannumeral 1) nor (\romannumeral 2), and the proof of Theorem 2.15 is complete.

\end{document}